\documentclass[a4paper,11pt]{article}
\usepackage{amsfonts,amsmath,amsthm}

\newtheorem{prop}{Proposition}
\newtheorem{thm}[prop]{Theorem}

\newtheorem{coroll}[prop]{Corollary}

\theoremstyle{remark}
\newtheorem{rmk}[prop]{Remark}

\numberwithin{equation}{section}

\renewcommand{\P}{\mathbb{P}}
\newcommand{\E}{\mathbb{E}}
\newcommand{\erre}{\mathbb{R}}
\newcommand{\cp}[2]{\langle#1,#2\rangle}
\newcommand{\bip}[2]{\big\langle#1,#2\big\rangle}
\newcommand{\ds}{\displaystyle}
\newcommand{\tr}{\mathop{\mathrm{Tr}}\nolimits}

\title{Well-posedness and invariant measures for HJM models with
  deterministic volatility and L\'evy noise}

\author{Carlo Marinelli\thanks{Institut f\"ur Angewandte Mathematik,
    Universit\"at Bonn, Wegelerstr. 6, D-53115 Bonn, Germany, and
    Centre de Recerca Matem\`atica, Bellaterra, Spain. E-mail:
    \texttt{marinelli@wiener.iam.uni-bonn.de}}}

\date{\normalsize First draft: December 29, 2006 \\
This version: May 12, 2008}

\begin{document}

\maketitle

\begin{abstract}
  We give sufficient conditions for existence, uniqueness and
  ergodicity of invariant measures for Musiela's stochastic partial
  differential equation with deterministic volatility and a Hilbert
  space valued driving L\'evy noise. Conditions for the absence of
  arbitrage and for the existence of mild solutions are also
  discussed.

\end{abstract}

\section{Introduction}
The aim of this work is to study some asymptotic properties of a
Heath-Jarrow-Morton model of the term structure of interest rates
driven by an infinite dimensional L\'evy noise. In particular,
denoting by $u(t,x)$, $t,\,x\geq 0$, the forward rate at time $t$
with maturity $t+x$, we shall be concerned with the following
stochastic partial differential equation (SPDE):
\begin{equation}
  \label{eq:spde}
du(t,x) = [u_x(t,x) + f(t,x)]\,dt + \cp{\sigma(t,x)}{dY_0(t)},
\end{equation}
where $Y_0$ is a L\'evy process taking values in a Hilbert space with
inner product $\cp{\cdot}{\cdot}$, $\sigma$ is a deterministic
volatility term, and $f$ is such that discounted prices of zero-coupon
bonds are local martingales. Precise assumptions will be stated below.
Note that in the special case where $Y_0$ is a Wiener process,
(\ref{eq:spde}) can be written in the more familiar form
\begin{equation}
  \label{eq:sacher}
du(t,x) = [u_x(t,x) + f(t,x)]\,dt
+ \sum_{k=1}^\infty \sigma^k(t,x)\,dw_k(t),
\end{equation}
where $w_k$ are real independent Wiener processes and $f$ satisfies
the well-known HJM drift condition (\cite{HJM}, \cite{filipo},
\cite{GM-review})
\begin{equation}
\label{eq:drift}
f(t,x) = \sum_{k=1}^\infty \sigma^k(t,x)\int_0^x\sigma^k(t,y)\,dy.
\end{equation}
Invariant measures and the asymptotic behavior of (\ref{eq:sacher}),
in the time-homogeneous case, have been studied by several authors
(see e.g. \cite{musiela}, \cite{vargiolu}, \cite{tehranchi}), allowing
also the volatility coefficient to depend on the forward rate itself.
Indeed it is widely accepted that mean reversion is a characteristic
property of the dynamics of interest rates, and it is supported by
empirical findings.

On the other hand, the literature on HJM models driven by L\'evy
processes has considerably grown in the last few years: let us just
cite, among others, \cite{ALM}, \cite{EbeJacRai}, and of course the
work from where it all started \cite{BDMKR} (where general random
measures are added to Brownian motion as driving noises). The
asymptotic behavior of such models, however, does not seem to have
been addressed. The present paper offers a first step in this
direction, in a simple model with deterministic volatility. The
setting is similar to that of \cite{vargiolu} (where the case of HJM
models driven by Wiener process was considered), but the choice of
state space is different, and the well-posedness of the model is
somewhat more complicated.

The paper is organized as follows: in section \ref{sec:muesli} we
derive sufficient conditions on $f$ ensuring that the bond market is
arbitrage free and write the SPDE (\ref{eq:spde}) as an abstract
evolution equation in a suitable Hilbert space of forward curves,
about which we discuss existence and uniqueness of mild and weak
solutions.
In section \ref{sec:conato} we discuss existence and uniqueness of
invariant measures, as well as the convergence in law of forward
curves as time goes to infinity.

After the first draft of this paper was completed, the author was
informed of some papers (most of which not yet published nor posted to
any standard preprint server) with some overlap with the present one.
In particular, drift conditions were derived in \cite{OezSchm} and in
\cite{JZ}. Our derivation, as well as the setup, is slightly
different, and we include it for the sake of completeness. Existence
and uniqueness of local mild solutions for HJM models driven by
finitely many independent L\'evy processes is discussed in
\cite{FilTap} (see also \cite{due}), where the volatility is also
allowed to be state dependent. In the present paper we consider a
possibly infinite dimensional driving L\'evy process, which is not
supposed to be the superposition of independent one-dimensional L\'evy
processes, and we look for global solutions. Moreover, in \cite{PZ07}
the authors obtain conditions for existence and uniqueness of global
mild solutions with an infinite dimensional L\'evy process and state
dependent volatility. However, their state space is a weighted $L_2$
space, which seems inappropriate for modeling purposes: in fact very
irregular forward curves cannot be excluded. We would also like to
point out that the choice of state space in \cite{FilTap} is the same
as in our paper, and for this reason the former paper only treats
local solutions. On the other hand, as already observed (see also
\cite{filipo} for an extensive discussion), our choice of state space
is better than that in \cite{PZ07}. Finally, the existence and
uniqueness statement of proposition \ref{prop:antipatica} can be
deduced by the result in \cite{Rus1} (see also \cite{Rus2} for a
result on ergodicity of HJM models with one-dimensional, possiblly
multiplicative noise, in spaces of weighted $L_2$ curves).

Let us conclude introducing some notation. Given two separable Hilbert
spaces $H$, $K$ we shall denote by $\mathcal{L}(H,K)$,
$\mathcal{L}_1(H,K)$ and $\mathcal{L}_2(H,K)$ the space of bounded
linear, trace-class, and Hilbert-Schmidt operators, respectively, from
$H$ to $K$.  $\mathcal{L}_1^+$ stands for the subset of
$\mathcal{L}_1$ of positive operators. We shall write
$\mathcal{L}_1(H)$ in place of $\mathcal{L}_1(H,H)$, and similarly for
the other spaces. Given a self-adjoint operator
$Q\in\mathcal{L}_1^+(H)$, we set $|x|_Q^2:=\cp{Qx}{x}$, $x\in H$. The
Hilbert-Schmidt norm is denoted by $|\cdot|_2$.
The characteristic function of a set $A$ is denoted by $\chi_A$, and
$\chi_r$ stands for the characteristic function of the set
$B_r:=\{x\in H:\;|x|\leq r\}$, where $H$ is a Hilbert space.
Given a continuously differentiable increasing function
$\alpha:\erre_+\to[1,\infty)$ such that $\alpha^{-1/3}\in
L_1(\erre_+)$, we define $L_{2,\alpha}^n:=L_{2,\alpha}^n(\erre_+)$ as
the space of distributions $\phi$ on $\erre_+$ such that
$\int_0^\infty |\phi^{(n)}(x)|^2\alpha(x)\,dx < \infty$.

\section{Musiela's SPDE with L\'evy noise}
\label{sec:muesli}
Throughout the paper $Y_0(t)$, $t\geq 0$, shall denote a L\'evy
process taking values in a (fixed) Hilbert space $K$, with generating
triplet $[b_0,R_0,m_0]$, i.e.
$$
\log\E e^{i\cp{y}{Y_0(1)}} =
i\cp{b_0}{y}
- \frac12 \cp{R_0y}{y}
+ \int_H (e^{i\cp{\xi}{y}}-1-i\cp{\xi}{y}\chi_1(\xi))\,m_0(d\xi),
$$
with $b_0\in K$, $R_0\in \mathcal{L}_1^+(K)$, and $m_0$ a
$\sigma$-finite measure on $\mathcal{B}(K)$, the Borel
$\sigma$-algebra of $K$, satisfying
$$
m_0(\{0\})=0, \quad
\int_K |\xi|^2 \,m_0(d\xi) < \infty
$$
(see e.g. \cite{GS-II} for details). This integrability assumption
serves two purposes: it ensures that $\E|Y_0(t)|^2$ is finite for all
$t\geq 0$, thus allowing to construct mild solutions of the SPDE
(\ref{eq:spde}) via an $L_2$ theory of stochastic integration, and it
allows to use Fubini's theorem to establish no-arbitrage sufficient
conditions. The assumption will turn out not to be a real restriction,
as the no-arbitrage condition essentially requires existence of
exponential moments of $m_0$ (see also \cite{ALM}).

\subsection{Drift condition}
\label{subsec:drift}
Let us consider, in the spirit of the original paper \cite{HJM} (see
also \cite{GM-review}), the following integral equation
\begin{equation}
\label{eq:kartoshka}
  u(t,x) = u(0,x) + \int_0^t b(s,x)\,ds 
                    + \int_0^t \cp{\sigma(s,x)}{dY_0(s)},
\end{equation}
where $b$ and $\sigma$ are random vector fields predictable in $s$ and
Borel measurable in $x$. In particular, they could depend on $u$
itself.

We shall give conditions under which the dynamics (\ref{eq:kartoshka})
is compatible with a no-arbitrage hypothesis, namely that the
corresponding discounted bond prices are local martingales. The
arguments of the proof consist of the L\'evy-It\^o decomposition in
Hilbert space (see e.g. \cite{Dett}) and the calculus for random
measures, following \cite{BDMKR}. We denote by $\hat{P}(t,\tau)$, $0\leq
t\leq \tau$, the discounted price of a zero-coupon bond expiring at
time $\tau$.
\begin{thm}
  Let $\Sigma(t,x)=-\int_0^x \sigma(t,y)\,dy$ and $E(t,x,\xi)=
  -\int_0^x \cp{\sigma(t,y)}{\xi}\,dy$ for all $x\geq 0$, $\xi \in K$.
  Assume that for all $t<\infty$
\begin{equation}
  \label{eq:h3}
  \int_0^t\int_0^\infty|b(s,x)|\,dx\,ds<\infty,
  \qquad
  \int_0^t\int_0^\infty|\sigma(s,x)|^2\,dx\,ds<\infty.
\end{equation}
Moreover, assume that for all $x\geq 0$ one has
\begin{equation}
\label{eq:locmg}
\begin{split}
\int_0^x b(t,y)\,dy =\;& u(t,x) - u(t,0)
+ \frac12 |\Sigma(t,x)|_{R_0}^2\\
& +\int_H \Big(e^{E(t,x,\xi)}-1-E(t,x,\xi)\chi_1(\xi)\Big)\,m_0(d\xi)  
\end{split}
\end{equation}
$d\P \times dt$-a.e. Then the discounted bond price process $t\mapsto
\hat{P}(t,\tau)$, $t\geq 0$, is a local martingale for all $\tau\geq
t$.
\end{thm}
\begin{proof}
The L\'evy process $Y_0$ admits the decomposition
\begin{equation}
  \label{eq:decomp}
  Y_0(t) = b_0t + W(t) + \int_0^t\int_{|\xi|<1} \xi\,\tilde{N}(ds,d\xi)
        + \int_0^t\int_{|\xi|\geq 1} \xi\,N(ds,d\xi),
\end{equation}
where $b_0 \in K$, $W$ is a $K$-valued Wiener process with covariance
operator $R_0$, $N$ is a Poisson measure on $K$ with compensator
$m_0$, and $\tilde{N}(ds,d\xi):=N(ds,d\xi)-ds\,m_0(d\xi)$.

Using the decomposition (\ref{eq:decomp}) one can write
\begin{eqnarray*}
  u(t,x) &=& u(0,x) + \int_0^t b(s,x)\,ds
                    + \int_0^t \cp{\sigma(s,x)}{dW(s)}\\
&&\qquad + \int_0^t \cp{\sigma(s,x)}{d\tilde{z}(s)}
         + \int_0^t \cp{\sigma(s,x)}{dz(s)},
\end{eqnarray*}
where
$$
\tilde{z}(t) = \int_0^t\int_{|\xi|< 1} \xi\,\tilde{N}(ds,d\xi),
\qquad
z(t) = \int_0^t\int_{|\xi|\geq 1} \xi\,N(ds,d\xi).
$$
However, let $(e_k)_{k\in\mathbb{N}}$ be a base of $K$ and set
$x^k=\cp{x}{e_k}$, $x\in K$. Then one has
\begin{eqnarray}
\int_0^t \cp{\sigma(s,x)}{d\tilde{z}(s)} &=&
\lim_{n\to\infty} \sum_{k=1}^n \int_0^t \sigma^k(s,x)\,d\tilde{z}^k(s) \nonumber\\
&=& \lim_{n\to\infty} \sum_{k=1}^n \int_0^t \int_{|\xi|<1}
            \sigma^k(s,x)\xi^k\,\tilde{N}(ds,d\xi) \nonumber\\
\label{eq:pi}
&=& \lim_{n\to\infty} \int_0^t \int_{|\xi|<1}
            \sum_{k=1}^n \sigma^k(s,x)\xi^k\,\tilde{N}(ds,d\xi),
\end{eqnarray}
where in the second line we have used the associativity of stochastic
integrals (i.e., using the notation of the stochastic calculus of
semimartingales and stochastic measures, $H\cdot(K\ast\mu) =
(HK)\ast\mu$ -- see e.g. \cite{JacShi}). Since
$\sum_{k=1}^n\sigma^k(s,x)\xi^k$ converges (for any fixed $x$) to
$\cp{\sigma(s,x)}{\xi}$ in $L_2(\Omega\times[0,t]\times
H,\mathcal{P},d\P\times ds\times dm_0)$ for all $t\geq 0$, with
$\mathcal{P}$ the predictable $\sigma$-field, one finally obtains
$$
\int_0^t \cp{\sigma(s,x)}{d\tilde{z}(s)} =
\int_0^t\int_{|\xi|< 1} \eta(s,x,\xi)\,\tilde{N}(ds,d\xi),
$$
where $\eta(s,x,\xi)=\cp{\sigma(s,x)}{\xi}$. Similar reasoning shows
that the same type of identity holds for integrals with respect to
$N$.
\par\noindent
For $x\geq 0$, let us denote by $p(t,x)$ the discounted prices at time
$t$ of a risk free zero-coupon bond expiring at time $t+x$.
By the definition of $p(t,x)$ and the equation for $u$ we obtain
\begin{eqnarray*}
  \log p(t,x) &=& -\int_0^t u(s,0)\,ds - \int_0^x u(t,y)\,dy \\
              &=& -\int_0^t u(s,0)\,ds - \int_0^x u(0,y)\,dy \\
              &&  -\int_0^x\int_0^t b(s,y)\,ds\,dy
                  -\int_0^x\int_0^t \cp{\sigma(s,y)}{dW(s)}\,dy \\
 && -\int_0^x\!\!\int_0^t\!\!\int_{|\xi|< 1}\eta(s,y,\xi)\,\tilde{N}(ds,d\xi)\,d
y
    -\!\int_0^x\!\!\int_0^t\!\!\int_{|\xi|\geq 1}\eta(s,y,\xi)\,N(ds,d\xi)\,dy \\
              &=& -\int_0^t u(s,0)\,ds - \int_0^x u(0,y)\,dy
              +\int_0^t B(s,x)\,ds + \int_0^t\cp{\Sigma(s,x)}{dW(s)} \\
   &&  +\int_0^t\int_{|\xi|< 1} E(s,x,\xi)\,\tilde{N}(ds,d\xi)
       +\int_0^t\int_{|\xi|\geq 1} E(s,x,\xi)\,N(ds,d\xi),
\end{eqnarray*}
where the third equality follows by Fubini's theorem (see \cite{DZ92} and 
\cite{lebedev}) and the definitions of $\Sigma$ and $E$, together with
$B(t,x):=-\int_0^t b(t,y)\,dy$.

Applying It\^o's formula, setting $\zeta(t,x)=\log p(t,x)$, one gets
\begin{eqnarray*}
  p(t,x) &=& e^{\zeta(t,x)} = e^{\zeta(0,x)} 
             + \int_0^t e^{\zeta(s-,x)} \Big(
                  -u(s,0)\,ds+B(s,x)\,ds + \cp{\Sigma(s,x)}{dW(s)}\Big) \\
         &&  + \frac12 \int_0^t e^{\zeta(s-,x)}\,d[\zeta,\zeta]^c(s,x)
             + \int_0^t\int_{|\xi|<1}%
                 (e^{\zeta(s-,x)+E(s,x,\xi)}-e^{\zeta(s-,x)})\,\tilde{N}(ds,d\xi
) \\
         &&  + \int_0^t\int_{|\xi|\geq 1}%
                 (e^{\zeta(s-,x)+E(s,x,\xi)}-e^{\zeta(s-,x)})\,N(ds,d\xi) \\
         &&  + \int_0^t\int_{|\xi|<1}%
                 (e^{\zeta(s-,x)+E(s,x,\xi)}-e^{\zeta(s-,x)}%
                      -E(s,x,\xi)e^{\zeta(s-,x)})\,m_0(d\xi)\,ds,
\end{eqnarray*}
or equivalently
\begin{eqnarray*}
\frac{dp(t,x)}{p(t-,x)} &=& \Big(-u(t,0)+B(t,x)
                              + \frac12|\Sigma(t,x)|_{R_0}^2\Big)dt
                              + \cp{\Sigma(t,x)}{dW(t)} \\
    && + \int_{|\xi|< 1}(e^{E(t,x,\xi)}-1)\,\tilde{N}(dt,d\xi) \\
    && + \int_{|\xi|\geq 1}(e^{E(t,x,\xi)}-1)\,N(dt,d\xi) \\
    && + \int_{|\xi|< 1}(e^{E(t,x,\xi)}-1-E(t,x,\xi))\,m_0(d\xi)\,dt.
\end{eqnarray*}
For $\tau\geq t$, by $\hat{P}(t,\tau)=p(t,\tau-t)$ it follows that
$d\hat{P}(t,\tau)=dp(t,\tau-t)-p_x(t,\tau-t)$, and 
$p(t,\tau-t)=e^{-\int_0^tu(s,0)\,ds}P(t,\tau-t)$ implies
$p_x(t,\tau-t)=e^{-\int_0^tu(s,0)\,ds}P_x(t,\tau-t)$.
Setting $P(t,\tau)=e^{-\int_0^\tau u(t,y)\,dy}$, one has
$P_x(t,\tau-t)=-P(t,\tau-t)u(t,\tau-t)$ and
$p_x(t,\tau-t)=-u(t,\tau-t)p(t,\tau-t)$, and finally
$$
d\hat{P}(t,\tau)=dp(t,\tau-t)+u(t,\tau-t)\hat{P}(t,\tau)\,dt.
$$
Together with the equation for $p(t,x)$, this implies
\begin{eqnarray*}
\frac{d\hat{P}(t,\tau)}{\hat{P}(t-,\tau)} &=&  \Big(-u(t,0) + u(t,\tau-t) + B(t,\tau-t)
                              + \frac12|\Sigma(t,\tau-t)|_{R_0}^2\Big)dt \\
    && + \cp{\Sigma(t,\tau-t)}{dW(t)}
            + \int_{|\xi|< 1}(e^{E(t,\tau-t,\xi)}-1)\,\tilde{N}(dt,d\xi) \\
    && + \int_{|\xi|\geq 1}(e^{E(t,\tau-t,\xi)}-1)\,N(dt,d\xi) \\
    && + \int_{|\xi|< 1}(e^{E(t,\tau-t,\xi)}-1-E(t,\tau-t,\xi))\,m_0(d\xi)\,dt.
\end{eqnarray*}
Lightening notation a bit, one can write
$$
\int_{|\xi|\geq 1} (e^E-1)\,dN =
\int_{|\xi|\geq 1} (e^E-1)\,d\tilde{N} 
+ \int_{|\xi|\geq 1} (e^E-1)\,dm_0\,dt,
$$
hence $\hat{P}(t,\tau)$ is a local martingale if
\begin{eqnarray*}
0 &=& -u(t,0) + u(t,\tau-t) + B(t,\tau-t)
+ \frac12|\Sigma(t,\tau-t)|_{R_0}^2 \\
&& \qquad + \int_{|\xi|\geq 1} (e^{E(t,\tau-t,\xi)}-1)\,m_0(d\xi)
+ \int_{|\xi|< 1} (e^{E(t,\tau-t,\xi)}-1-E(t,\tau-t,\xi))\,m_0(d\xi) \\
&=& -u(t,0) + u(t,\tau-t) + B(t,\tau-t)
+ \frac12|\Sigma(t,\tau-t)|_{R_0}^2 \\
&& \qquad + \int_K \big(e^{E(t,\tau-t,\xi)}-1-E(t,\tau-t,\xi)\chi_1(\xi)\big)\,m_0(d\xi),
\end{eqnarray*}
and the theorem is proved.
\end{proof}

\begin{rmk}
  The above theorem implies a ``drift condition'' that generalizes the
  HJM condition (\ref{eq:drift}). In particular, assume that $Y(1)$
  admits exponential moments, or equivalently that $\int_H
  e^{\cp{z}{\xi}}\,m_0(d\xi)<\infty$ for all $z\in H$, and define the
  function $\psi(z)=\log\,\E e^{\cp{z}{Y(1)}}$. Then taking into
  account (\ref{eq:spde}) and (\ref{eq:kartoshka}), (\ref{eq:locmg})
  implies the following relation between the drift and the volatility
  functions:
  \begin{equation}
    \label{eq:jj}
    \int_0^\cdot f(t,y)\,dy = \psi\Big(-\int_0^\cdot \sigma(t,y)\,dy\Big)
  \end{equation}
  for all $t\in[0,T]$. Unfortunately this identity is ``implicit'',
  and only under further assumptions can it be made more explicit (see
  below).
\end{rmk}

\subsection{Abstract setting and well-posedness}
We shall rewrite the SPDE (\ref{eq:spde}) as an abstract stochastic
differential equation in the space $H=L_{2,\alpha}^1$. The space $H$
endowed with the inner product
$$
\cp{\phi}{\psi} = \int_{\erre_+} \phi'(x)\psi'(x)\,\alpha(x)\,dx
+ \lim_{x\to\infty} \phi(x)\psi(x)
$$
is a separable Hilbert space. This choice of state space is standard
and is apparently due to Filipovi\'c \cite{filipo}. Nonetheless, other
authors have studied related SPDEs in different function spaces, e.g.
in weighted $L_2$ spaces, weighted Sobolev spaces, or fractional
Sobolev spaces (see \cite{GM-review}, \cite{vargiolu},
\cite{EkeTaf} respectively).

Let us define on $H$ the operator $A:f\mapsto f'$, with domain
$D(A)=L^1_{2,\alpha}\cap L^2_{2,\alpha}$, which generates the
semigroup of right shifts $[e^{tA}\phi](x):=\phi(x+t)$, $t\geq 0$.
Musiela's SPDE (\ref{eq:spde}) can be written in abstract form as
\begin{equation}
\label{eq:muso}
du(t) = (Au(t) + f(t))\,dt + B(t)\,dY_0(t),
\end{equation}
where $f(t) \equiv f(t,\cdot)$ and $B(t)\in\mathcal{L}(K,H)$ is
defined by $[B(t)u](\cdot)=\cp{\sigma(t,\cdot)}{u}_K$, with suitable
regularity assumptions on $\sigma$.

Several papers deal with the solution of this type of equations in the
time-independent case with $f(t)\equiv 0$, $B(t)\equiv B$. Here we
limit ourselves to mention \cite{choj}, which is probably the first
paper considering weak solutions (in the sense of PDEs), and
\cite{FR}, where an analytic approach is used to solve, even in the
strong sense, equations of the type (\ref{eq:museto}), possibly in a
larger space, allowing the characteristic function of $Y_0(1)$ to be
only continuous instead of Sazonov continuous (see also
\cite{LR-pseudo}, \cite{LR-heat}, \cite{RW-harnack}).

Our goal in this subsection is less ambitious, namely we shall only
prove that the formal solution of (\ref{eq:muso}) given by the
variation of constants formula
\begin{equation}
\label{eq:formal}
u(t) = e^{tA}u_0 + \int_0^t e^{(t-s)A}f(s)\,ds + \int_0^t e^{(t-s)A}B(s)\,dY_0(s)
\end{equation}
is a well defined process and provides the unique weak solution to
(\ref{eq:muso}). In particular, this implies that the L\'evy-based
model (\ref{eq:spde}) for the evolution of forward curves is well
posed under appropriate assumptions on $\sigma$.

In analogy to \cite{DP-K}, we shall prove that (\ref{eq:formal})
belongs to the space $\mathcal{H}_2(T)$ of mean square continuous
process on $[0,T]$ with values in $H$ adapted to the filtration
generated by $Y$. We endow $\mathcal{H}_2(T)$ with the norm defined by
$$
\|F\|^2_2 := \sup_{t\in[0,T]} \E|F(t)|_H^2.
$$

\begin{prop}\label{prop:tiepida}
  Assume that
  $$
  \int_0^T (|f(t)|_H + |B(t)|^2_2)\,dt < \infty.
  $$
  Then $u$ defined as in (\ref{eq:formal}) belongs to
  $\mathcal{H}_2(T)$.
\end{prop}
\begin{proof}
Adaptedness is immediate by definition. Using the L\'evy-It\^o
decomposition in the form
$$
Y_0(t) = at + W(t) + \int_0^t\int_K \xi \tilde{N}(ds,d\xi),
$$
with $a=b_0 + \int_{|\xi|\geq 1} \xi\,m_0(d\xi)$, and taking into
account Propositions 2.2, 2.3 in \cite{DP-K}, it is enough to consider
the case when $Y_0(t)$ has no drift and no Brownian component.
In particular one has
\begin{equation*}
\E\bigg|\int_0^t e^{(t-s)A}B(s)\int_K \xi\tilde{N}(ds,d\xi)\bigg|^2
\leq N^2
\int_0^T|B(t)|_2^2\int_K
|\xi|^2\,m_0(d\xi)\,dt < \infty,
\end{equation*}
with $N=\sup_{t\in[0,T]} |e^{tA}|$.  Moreover $|\int_0^t
e^{(t-s)A}f(s)\,ds| \leq N\int_0^T|f(s)|\,ds<\infty$, hence
$\|u\|_2<\infty$.

Let us now prove that $t\mapsto \E|u(t)|^2$ is continuous. Setting
$Y_A(t):=\int_0^t e^{(t-s)A}B(r)\,dY_0(r)$, it is enough to prove that
$t\mapsto \E|Y_A(t)|^2$ is continuous. For $0\leq s\leq t\leq T$, write
$$
Y_A(t)-Y_A(s) = \int_0^s \big(e^{(t-r)A}-e^{(s-r)A}\big)B(r)\,dY_0(r)
+ \int_s^t e^{(t-r)A}B(r)\,dY_0(r),
$$
where the two terms on the right-hand side are uncorrelated.
Since
\begin{equation*}
\begin{split}
& \E\left|\int_0^s \big(e^{(t-r)A}-e^{(s-r)A}\big)B(r)\,dY_0(r)\right|^2
\\
& \phantom{\E\left| \int_s^t e^{(t-r)A}B(r)\,dY_0(r) \right|^2}
\leq
\int_0^s \big|e^{(t-r)A}-e^{(s-r)A}\big|^2|B(r)|_2^2\,dr
     \int_K |\xi|^2\,m_0(d\xi) \to 0,\\
& \E\left| \int_s^t e^{(t-r)A}B(r)\,dY_0(r) \right|^2
\leq
\int_s^t \big|e^{(t-r)A}\big|^2\,|B(r)|_2^2\,dr \int_K |\xi|^2\,m_0(d\xi) \to 0
\end{split}
\end{equation*}
as $s \uparrow t$, the result follows. The case $0 \leq t \leq s \leq
T$ and $s \downarrow t$ is completely analogous, hence omitted.
\end{proof}

The following proposition shows that existence in the mild sense for
equation (\ref{eq:muso}) implies existence and uniqueness in the weak
sense. This fact was essentially proved by A.~Chojnowska-Michalik in
\cite{choj} in the time-independent case. Here we give a more direct
proof that closely follows \cite{DZ92}.
\begin{prop}
  Equation (\ref{eq:muso}) has a unique weak solution given by
  (\ref{eq:formal}).
\end{prop}
\begin{proof}
Let us define the additive process $Y(t)=f(t)+\int_0^tB(s)\,dY_0(s)$
and write (\ref{eq:muso}) in the form
$$
du(t) = Au(t)\,dt + dY(t), \qquad u(0)=u_0.
$$
It is enough to consider the case $u_0=0$, the extension to the
general case being immediate. We need to prove that, for $v\in
D(A^*)$,
$$
\bip{\int_0^t e^{(t-s)A}\,dY(s)}{v} =
\int_0^t \bip{\int_0^s e^{(s-r)A}\,dY(r)}{A^*v}\,ds
+ \cp{Y(t)}{v},
$$
or equivalently
$$
\bip{\int_0^t (e^{(t-s)A}-I)\,dY(s)}{v} =
\int_0^t \bip{\int_0^s e^{(s-r)A}\,dY(r)}{A^*v}\,ds.
$$
In fact, the right-hand side of the previous expression can be
equivalently written as
\[
\int_0^t\int_0^s \big\langle
e^{(s-r)A^*}A^*v,dY(r)\big\rangle\,ds,
\]
hence, using the identity
\[
\int_0^{t-s} e^{rA^*}A^*v\,dr = \big(e^{(t-s)A^*}-I\big)v,
\]
Fubini's theorem yields the desider conclusion. Uniqueness will follow
if we prove that a weak solution is a mild solution. One immediately
recognize that the proof of Lemma 5.5 in \cite{DZ92}, repeated word by
word, yields the identity
$$
\cp{w(t)}{\phi(t)} = \int_0^t \cp{w(s)}{\phi'(s)+A^*\phi(s)}\,ds
+ \int_0^t \cp{\phi(s)}{dY(s)}, \quad t\in[0,T],
$$
where $w$ is a weak solution of (\ref{eq:museto}) and $\phi\in
C^1([0,T],D(A^*))$. Taking $\phi(s)=e^{(t-s)A^*}\phi_0$, $\phi_0\in
D(A^*)$, implies
$$
\cp{w(t)}{\phi_0} = \bip{\int_0^t e^{(t-s)A}dY(s)}{\phi_0},
$$
hence $u(t)=w(t)$ because $D(A^*) \subset H$ densely.
\end{proof}

It is clear that in order to obtain results on the well-posedness of
the Musiela's SPDE (\ref{eq:spde}) it is necessary to establish
conditions on $\sigma$ and $Y_0$ such that the hypotheses of
Proposition \ref{prop:tiepida} are satisfied. The following sufficient
conditions are obviously the most general ones, but also the hardest to
verify.
\begin{prop} Let $(e_k)_{k\in\mathbb{N}}$ be a basis of $K$, and define
  $\sigma^k(t,x)=\cp{\sigma(t,x)}{e_k}_K$. Assume that
  \begin{equation}
  \label{eq:serrano}
  \int_0^T \sum_{k=1}^\infty \int_0^\infty |\sigma^k_x(t,x)|^2\alpha(x)\,dx\,dt
  < \infty
  \end{equation}
  and that there exists $f:[0,T] \to H$ satisfying (\ref{eq:jj}) such
  that $\int_0^T |f(t)|\,dt<\infty$. Then there exists the mild
  solution of (\ref{eq:muso}).
\end{prop}
\begin{proof}
  One has $|B(t)|_2^2=\sum_{k=1}^\infty |B(t)e_k|^2$, and
  $B(t)e_k=\sigma^k(t,\cdot)$, hence
  $$
  |B(t)e_k|^2=\int_0^\infty |\sigma^k_x(t,x)|^2\alpha(x)\,dx.
  $$
  Therefore (\ref{eq:serrano}) implies $\int_0^T
  |B(t)|_2^2\,dt<\infty$, and the result follows by proposition
  \ref{prop:tiepida}.
\end{proof}

Much simpler conditions can be stated if the driving noise $Y_0$ is
finite dimensional.
\begin{prop}
  Assume that $Y_0$ is a $\erre^n$-valued L\'evy process such that
  $\int_{\erre^n}e^{\cp{\xi}{z}}\,m_0(d\xi)<\infty$ for all $z\in
  B_r$, for some $r>0$. If $\int_0^x\sigma(t,y)\,dy \in B_r$ for all
  $(t,x)\in[0,T]\times\erre_+$, and
  \begin{equation}
  \label{eq:parma}
  \int_0^T\Big(\int_0^\infty |\sigma^k_x(t,x)|^2\alpha(x)\,dx\Big)^2\,dt
  < \infty
  \end{equation}
  for all $k=1,\ldots,n$, then (\ref{eq:jj}) reduces to
  \begin{equation}
    \label{eq:deda}
    f(t,x) = -\Big\langle \sigma(t,x),
           D\psi\Big(-\int_0^x\sigma(t,y)\,dy\Big) \Big\rangle,
    \qquad (t,x)\in [0,T]\times\erre_+,
  \end{equation}
  and (\ref{eq:muso}) admits a mild solution.
\end{prop}
\begin{proof}
  It is enough to check that the hypotheses of proposition
  \ref{prop:tiepida} are satisfied. In particular, similarly as
  before, one has, using Jensen's inequality,
  \begin{align*}
  \int_0^T |B(t)|_2^2\,dt &= \int_0^T \sum_{k=1}^n \int_0^\infty
        |\sigma^k_x(t,x)|^2\alpha(x)\,dx\,dt\\
  &=\sum_{k=1}^n  \int_0^T |\sigma^k(t,\cdot)|_H^2\,dt
     \leq \sum_{k=1}^n \Big( T\int_0^T |\sigma^k(t,\cdot)|_H^4\,dt \Big)^{1/2}
     < \infty
  \end{align*}
  where the last inequality is immediate by (\ref{eq:parma}).
  Moreover, in (\ref{eq:deda}) the quantity $D\psi(x)$ is well defined
  for $|x|\leq r$ because $\psi\in C^\infty(B_r)$. This
  implies that there exists a positive constant $N$ such that
  $|\psi_{x^i}(z)|<N$, $|\psi_{x^ix^j}(z)|<N$ for all $z\in B_r$. We
  have
  \begin{eqnarray*}
  f_x(t,x) &=& -\Big\langle \sigma_x(t,x),
         D\psi\Big(-\int_0^x\sigma(t,y)\,dy\Big) \Big\rangle\\
           && + \Big\langle D^2\psi\Big(-\int_0^x\sigma(t,y)\,dy\Big)
                 \sigma(t,x),\sigma(t,x) \Big\rangle,
  \end{eqnarray*}
  hence $\int_0^T |f(t)|_H\,dt < \infty$ if $\int_0^T\int_0^\infty
  [\sigma^i(t,x)\sigma^j(t,x)]^2 \alpha(x)\,dx\,dt < \infty$ for all
  $i$, $j\leq n$.  The condition $\int_0^x\sigma(t,y)\,dy<r$ for all
  $x\geq 0$ implies that $\lim_{y\to\infty}|\sigma(t,y)|=0$ for all
  $t\in[0,T]$, hence, using Cauchy-Schwarz' inequality,
  \begin{eqnarray*}
  \lefteqn{\int_0^T\!\int_0^\infty [\sigma^i(t,x)\sigma^j(t,x)]^2 \alpha(x)\,dx\,dt}\\
  && \leq \Big(\int_0^T\!\int_0^\infty \sigma^i(t,x)^4\alpha(x)\,dx\,dt\Big)^{1/2}
  \Big(\int_0^T\!\int_0^\infty \sigma^j(t,x)^4\alpha(x)\,dx\,dt\Big)^{1/2}\\
  && \leq N(\alpha) \Big(\int_0^T |\sigma^i(t,\cdot)|_H^4\,dt\Big)^{1/2}
                    \Big(\int_0^T |\sigma^j(t,\cdot)|_H^4\,dt\Big)^{1/2}
  \end{eqnarray*}
  where the second inequality follows by (5.8) in \cite{filipo}.
\end{proof}

An analogous expression could be obtained for a general Hilbert space
valued noise $Y$, if one can guarantee that $\psi$ is Fr\'echet
differentiable. In the next proposition we shall identify the
Fr\'echet derivative $D\psi(x) \in \mathcal{L}(K,\erre)$ with its
Riesz representative vector in $K$.
\begin{prop}
  Assume that $\int_K e^{\cp{\xi}{z}}\,m_0(d\xi)<\infty$ for all $z\in
  B_r$, $\psi\in C^2_b(B_r)$, for some $r>0$, and
    \begin{equation}
  \label{eq:chorizo}
  \int_0^T \Big(\sum_{k=1}^\infty \int_0^\infty
        |\sigma^k_x(t,x)|^2\alpha(x)\,dx\Big)^2 dt < \infty
  \end{equation}
  If $\int_0^x\sigma(t,y)\,dy \in B_r$ for all
  $(t,x)\in[0,T]\times\erre_+$, then (\ref{eq:jj}) reduces to
  (\ref{eq:deda}), where $\cp{\cdot}{\cdot}$ is the inner product of
  $K$, and equation (\ref{eq:muso}) admits a mild solution.
\end{prop}
\begin{proof}
  As seen before, we have
  $$
  |B(t)|_2^2 = \sum_{k=1}^\infty \int_0^\infty
        |\sigma^k_x(t,x)|^2\alpha(x)\,dx,
  $$
  hence, by (\ref{eq:chorizo}), $\int_0^T
  |B(t)|_2^2\,dt\leq\big(T\int_0^T|B(t)_2^4\,dt\big)^{1/2}<\infty$.
  Let us prove that $\int_0^T |f(t)|\,dt<\infty$. The same expression
  for $f_x(t,x)$ as in the proof of the previous proposition holds,
  mutatis mutandis.  Therefore, in view of (\ref{eq:chorizo}) and
  $\psi\in C^2_b(B_r)$, it is enough to show that
  $\int_0^T\!\int_0^\infty |\sigma(t,x)|^4\alpha(x)\,dx\,dt<\infty$.
  Let $(e_k)$ be a basis of $K$, and set
  $\sigma_n(t,x)=\sum_{k=1}^n\sigma^k(t,x)e_k$. Let
  $\phi^{(\varepsilon)}_n$ a smooth approximation of $x\mapsto |x|$ in
  $\erre^n$ such that $|D\phi^{(\varepsilon)}_n| \leq 1$, then we have
  as follows by (5.8) in \cite{filipo}
  \begin{align*}
  \int_0^\infty \phi^{(\varepsilon)}_n(\sigma_n(t,x))^4\alpha(x)\,dx
  &\leq N \big|\phi^{(\varepsilon)}_n(\sigma_n(t,\cdot))\big|_H^4\\
  &=N\Big( \int_0^\infty \big|
          D_x\phi^{(\varepsilon)}_n(\sigma_n(t,x))\big|^2\alpha(x)\,dx\Big)^2\\
  &\leq N \Big( \int_0^\infty \big| D_x\sigma_n(t,x)\big|^2
                  \alpha(x)\,dx \Big)^2 \\
  &\leq N \Big( \sum_{k=1}^\infty \int_0^\infty
                  |\sigma^k_x(t,x)|^2\alpha(x)\,dx \Big)^2
  \end{align*}
  with $N=N(\alpha)$, thus by (\ref{eq:chorizo})
  \begin{equation}
      \label{eq:eppa}
  \int_0^T\!\!\int_0^\infty \phi^{(\varepsilon)}_n(\sigma_n(t,x))^4
              \alpha(x)\,dx\,dt \leq N
  \int_0^T\!\!\Big(\sum_{k=1}^\infty \int_0^\infty
        |\sigma^k_x(t,x)|^2\alpha(x)\,dx\Big)^2 dt < \infty.
  \end{equation}
  Since the bound in (\ref{eq:eppa}) does not depend on $\varepsilon$
  nor on $n$, passing to the limit as $\varepsilon \to 0$ we get
  $\int_0^T\!\int_0^\infty |\sigma_n(t,x)|^4\alpha(x)\,dx < \infty$, and letting
  $n$ tend to infinity we finally get $\int_0^T\!\int_0^\infty
  |\sigma(t,x)|^4\alpha(x)\,dx < \infty$.
\end{proof}


\section{Invariant measures and asymptotic behavior}
\label{sec:conato}
In this section we assume $\sigma(t,\cdot)\equiv \sigma(\cdot)$, thus
also $B(t)\equiv B$, $f(t)\equiv f$. In view of the no-arbitrage
considerations in the previous section, we also assume that $m_0$
admits exponential moments, and $B\in\mathcal{L}(K,H_0)$, where
$H_0:=\{g\in H:\;g(\infty)=0\}$, hence also $f\in H_0$. In order for
the following results to hold, it is not necessary to assume that $f$
is such that no-arbitrage is verified, even though this is of course
the situation we are interested in.

Let us rewrite (\ref{eq:muso}), for convenience of notation, in the
more compact form
\begin{equation}
\label{eq:museto}
du(t) = Au(t)\,dt + dY(t),
\end{equation}
where $Y(t):=ft + BY_0(t)$. Then one can easily prove that $Y$ is
a $H$-valued L\'evy process with triplet $[b,R,m]$, where
\begin{eqnarray*}
b &=& f + Bb_0 + \int_K B\xi\big(\chi_1(B\xi)-\chi_1(\xi)\big)\,m_0(d\xi)\\
R &=& BR_0B^* \\
m(d\xi) &=& m_0(B^{-1}d\xi)
\end{eqnarray*}

The following proposition gives a simple sufficient condition for the
existence and uniqueness of an invariant measure for an HJM model with
deterministic volatility and Hilbert space valued L\'evy noise. The
only real requirement is that the state space $L_{2,\alpha}^1$ is
chosen with an exponentially growing weight $\alpha$.

\begin{prop}     \label{prop:antipatica}
  Assume that $\alpha(x)=e^{\alpha x}$, $\alpha>0$, and the forward
  curve at time zero is deterministic. Then there exists a unique
  invariant measure for (\ref{eq:museto}) to which the law of $u(t)$
  weakly converges as $t\to\infty$.
\end{prop}
\begin{proof}
  Writing equation (\ref{eq:museto}) in mild form, recalling that the
  range of $B$ is contained in $H_0$, one recognizes that
  $u(t,\infty)=u_0(\infty)$ for all $t\geq 0$ (``long rates never
  fall''). Considering the isomorphism $H=H_0\oplus\erre$,
  (\ref{eq:museto}) is equivalent to the system
  \begin{equation}
  \label{eq:faxoi}
  \begin{cases}
    d\bar u(t) = A\bar u(t)\,dt + dY(t)\\
    u(t,\infty) = \ell,
  \end{cases}
  \end{equation}
  where $\bar u$ is the projection of $u$ on $H_0$, $A$ still denotes
  the restriction of $A$ to $H_0$, and $\ell\in\erre$.  Let us show
  that $e^{tA}$ is exponentially stable on $H_0$:
  $$
  |e^{tA}\phi|^2_{H_0} = \int_0^\infty \phi'(x+t)^2\alpha(x)\,dx
  \leq e^{-\alpha t} \int_0^\infty \phi'(x)^2\alpha(x)\,dx
  = e^{-\alpha t} |\phi|_{H_0}^2,
  $$
  i.e. $|e^{tA}|\leq e^{-t\alpha/2}$. The obvious inequality
  $x^2\geq\log(1+x)$, $x\geq 1$, and the assumption $\int_K
  |\xi|^2\,m_0(d\xi)<\infty$ imply that
  \begin{align*}
  \int_{|\xi|\geq 1} \log(1+|\xi|)\,m(d\xi) &\leq
  \int_{|\xi|\geq 1} |\xi|^2\,m_0(B^{-1}d\xi)\\
  &\leq |B|^2 \int_K |\zeta|^2\,m_0(d\zeta) < \infty.
  \end{align*}
  Therefore theorem 6.7 of \cite{choj} yields the existence of an
  invariant measure $\bar\mu$ on $H_0$ for the first equation of
  (\ref{eq:faxoi}), hence $\mu=\bar\mu\otimes\delta_\ell$ is an
  invariant measure for (\ref{eq:museto}) on $H$. Since $|e^{tA}|\leq
  e^{-t\alpha/2}\to 0$ as $t\to\infty$ (i.e. $e^{tA}$ is stable),
  proposition 6.1 of \cite{choj} (or theorem 3.1 of \cite{FR}) imply
  that $\bar\mu$ is infinitely divisible with triplet
  $[b_\infty,R_\infty,m_\infty]$, where
  \begin{align}
    \label{eq:binf}
    b_\infty &= \lim_{t\to\infty}\Big[ \int_0^t e^{sA}b\,ds 
       + \int_0^t\int_{H_0} e^{sA}\xi 
            \big(\chi_1(e^{sA}\xi)-\chi_1(\xi)\big)\,m(d\xi)\,ds\Big]\\
    R_\infty &= \int_0^\infty e^{sA} R e^{sA^*}\,ds\\
    m_\infty(d\xi) &= \int_0^\infty m((e^{sA})^{-1}d\xi)\,ds,
    \quad m_\infty(\{0\})=0, \label{eq:minf}
  \end{align}
  are all well defined thanks to the stability properties of $e^{tA}$.
  In particular $\bar\mu$ is unique. Finally, by proposition 6.1 of
  \cite{choj} we have that $\bar\mu$ coincides with the law of the
  random variable $\int_0^\infty e^{sA}\,dY(s)$, and lemma 3.1 of
  \cite{choj} allows to conclude that the law of $\bar u(t)$ weakly
  converges to $\bar\mu$.
\end{proof}
\begin{rmk}
  The decomposition $H=H_0\oplus\erre$ was already used in
  \cite{tehranchi}, but essentially the same ``trick'' already
  appeared, perhaps less explictly, in \cite{vargiolu}. The fast
  growth at infinity of the weight $\alpha$ was needed in
  \cite{tehranchi} as well, where it is assumed that
  $\alpha:=\inf_{x\geq 0} \alpha'(x)/\alpha(x)>0$. In fact, Gronwall's
  lemma immediately yields that this condition implies
  $\alpha(x)\geq\alpha(0)e^{\alpha x}$.
\end{rmk}

The choice of the weight function $\alpha$, as just seen, determines
the stability properties of the semigroup $e^{tA}$. For a generic
choice of $\alpha$ we cannot guarantee exponential stability of
$e^{tA}$ in $H_0$, but we still have stability, in the sense that
$|e^{tA}g|_{H_0}\to 0$ as $t\to\infty$ for any $g\in H_0$. However, in
order to obtain existence of an invariant measure for
(\ref{eq:museto}), the conditions to verify become quite difficult, in
general. In particular the following characterization holds, the proof
of which follows \cite{choj} or \cite{FR}.
\begin{prop}
  Assume that the forward curve at time zero is deterministic. The
  following conditions are sufficient and necessary for the existence
  of a (unique) invariant measure $\mu$ for (\ref{eq:museto}):
  \begin{itemize}
  \item[\emph{(i)}] $\ds\sup_{t\geq 0} \tr \int_0^t e^{sA}BR_0B^*e^{sA^*}\,ds$;
  \item[\emph{(ii)}] $\ds\int_0^\infty\!\int_{H_0} (|e^{tA}x|^2\wedge 1)\,m(dx)\,dt<\infty$;
  \item[\emph{(iii)}] the limit in (\ref{eq:binf}) exists.
  \end{itemize}
  Moreover, $\mu=\bar\mu \otimes \delta_\ell$, where $\bar\mu$ is
  infinitely divisible with triplet $[b_\infty,R_\infty,m_\infty]$
  given by (\ref{eq:binf})-(\ref{eq:minf}). Finally, the law of $u(t)$
  weakly converges to $\bar\mu$ at $t\to\infty$.
\end{prop}
\begin{proof}
  The semigroup $e^{tA}$ is stable on $H_0$ because
  \begin{equation}
    \label{eq:contra}
    |e^{tA}g|^2_{H_0} = \int_t^\infty g'(x)^2 \alpha(x-t)\,dx
    \leq \int_t^\infty g'(x)^2 \alpha(x)\,dx
    \stackrel{t\to\infty}{\longrightarrow} 0,
  \end{equation}
  where the inequality follows by monotonicity of $\alpha$ and the
  limit is zero because the integrand is in $L_1(\erre_+)$.  Therefore
  theorem 6.4 of \cite{choj} (or theorem 3.1 of \cite{FR}) implies
  that the infinitely divisible measure $\bar\mu$ on $H_0$ with
  triplet $[b_\infty,R_\infty,m_\infty]$ is invariant for the first
  equation of (\ref{eq:faxoi}). The proof is then completed exactly as
  in the previous proposition.
\end{proof}

\begin{coroll}
\label{cor:boh}
  Assume that the forward curve at time zero is deterministic and that
  there exists a function $\phi\in L_1(\erre_+)\cap L_2(\erre_+)$ such
  that $|e^{tA}x|_{H_0} \leq \phi(t)|x|_{H_0}$. Then $\mu$ as defined
  in the previous proposition is the unique invariant measure of
  (\ref{eq:museto}) and it is ergodic.
\end{coroll}
\begin{proof}
  Since $e^{tA}$ is a stable semigroup on $H_0$, it is enough to
  verify hypotheses (i)--(iii) of the last proposition. We have
\begin{eqnarray*}
\int_0^\infty\!\!\int_{H_0} (|e^{tA}x|^2 \wedge 1)\, m(dx)\,dt &\leq&
\int_0^\infty\!\!\int_{H_0} \phi(t)^2|x|^2\, m(dx)\,dt\\
&\leq& \int_0^\infty \phi(t)^2\,dt \;
       |B|^{2}\!\int_K |x|^{2}\, m_0(dx) < \infty,
\end{eqnarray*}
because $\int_K |x|^2\,m_0(dx)<\infty$. Since $m_0$ admits exponential
moments, then $\int_K |x|\,m_0(dx)<\infty$, and
$$
\int_0^\infty\!\int_{H_0} \big|e^{tA}x
   \big(\chi_1(e^{tA}x)-\chi_1(x)\big)\big|\,m(dx)\,dt \leq
2\int_0^\infty \phi(t)\,dt |B| \int_K |x|m_0(dx) < \infty.
$$
Let us now prove that $\lim_{t\to\infty} \int_0^t
e^{sA}b\,ds$ exists in $H_0$: we have
$$
\bar b(x) := \lim_{t\to\infty}\Big[\int_0^t e^{sA}b\,ds\Big](x)
= \int_x^\infty b(s)\,ds, \qquad x\geq 0,
$$
thus $\bar{b}'(x)=-b(x)$ and
$$
\int_0^\infty |e^{sA}b|\,ds \leq
|b|_{H_0} \int_0^\infty \phi(s)\,ds < \infty,
$$
i.e. $b_\infty\in H_0$. Similarly we have
$$
\tr R_\infty = \sup_{t\geq 0} \int_0^t e^{sA}BR_0B^*e^{sA^*}\,ds
\leq |B|^2\tr R_0 \int_0^\infty \phi(s)^2\,ds < \infty,
$$
i.e. $R_\infty$ is well defined. 
\end{proof}
\begin{rmk}
  The semigroup $e^{sA}$ is a contraction semigroup for any choice of
  $\alpha$, e.g. by (\ref{eq:contra}), hence if $\phi$ as in the
  previous corollary exists, then one can always choose $|\phi(t)|\leq
  1$, thus $\phi\in L_1(\erre_+)$ also implies $\phi\in L_2(\erre_+)$.
  A possible choice of $\phi$ (although very rough) could be
  $$
  \phi(t) = \sup_{x\geq 0}
  \Big(\frac{\alpha(x)}{\alpha(x+t)}\Big)^{1/2},
  $$
  provided that $\int_0^\infty \sup_{x\geq 0}
  \big(\alpha(x)/\alpha(x+t)\big)^{1/2}\,dt<\infty$.
\end{rmk}

We have seen in the previous proposition that one of the necessary and
sufficient conditions for the existence of an invariant measure is
that $b_\infty$ exists. In particular, if $m_0$ is symmetric, the
problem reduces to proving that $\bar{b}:=\int_0^\infty e^{sA}b\,ds$
is a well-defined element of $H_0$. In fact, if $m_0$ is symmetric
then $m$ is symmetric as well, and the second summand on the right hand
side of (\ref{eq:binf}) is zero for all $t\geq 0$.  It is thus natural to
look for conditions on $\alpha$ such that the norm of $\bar{b}$ can be
bounded in terms of the norm of $b$. This is indeed possible, and one can give a sharp condition, namely (\ref{eq:mucken}) below is necessary and
sufficient for $|\bar{b}| \leq N |b|$ to hold.
\begin{prop}
Assume that 
\begin{equation}
  \label{eq:mucken}
  \sup_{x\geq 0} \int_0^x \alpha(y)\,dy
                 \int_x^\infty \frac1{\alpha(y)}\,dy < \infty.  
\end{equation}
Then $\int_0^\infty e^{sA}b\,ds$ exists in $H_0$.
\end{prop}
\begin{proof}
  As in the proof of corollary \ref{cor:boh}, we only need to prove
  that $b^2\alpha\in L_1(\erre_+)$. Let $\nu$ be a nonnegative Borel
  measure on $\erre_+$. By a result of Muckenhoupt \cite{mucken}, we
  have that the following weighted Hardy inequality holds for all
  measurable functions $f$
  \begin{equation*}
    \int_0^\infty \Big|\int_0^x f(y)\,dy\Big|^2 \,\nu(dx) \leq 
    N \int_0^\infty f(x)^2\,\nu(dx),
  \end{equation*}
  where $N$ is a positive constant that depends only on $\alpha$,
  if and only if
  $$
  \sup_{r\geq 0} [\nu([r,+\infty))]^{1/2} 
             \Big[\int_0^r \Big(\frac{d\nu}{dx}\Big)^{-1}\,dx\Big]^{1/2}
  < \infty.
  $$
  By the change of variable $x\mapsto x^{-1}$, choosing
  $\nu(dx)=\alpha(x)\,dx$, we obtain that (\ref{eq:mucken}) is
  necessary and sufficient for
  \begin{equation*}
    \int_0^\infty \Big|\int_x^\infty b'(y)\,dy\Big|^2 \alpha(x)\,dx \leq 
    N \int_0^\infty b'(x)^2 \alpha(x)\,dx.
  \end{equation*}
  Since $b\in H_0$, we have that $|\int_x^\infty
  b'(y)\,dy|=|b(x)|$ and that the right-hand side of the previous
  inequality is finite.
\end{proof}

\section{Conclusions}
We have considered an equation of HJM type driven by a L\'evy process
taking values in a Hilbert space, obtaining sufficient conditions for
the absence of arbitrage. Assuming that the volatility operator is
deterministic, and using Musiela's parametrization, one obtains a
stochastic evolution equation of Ornstein-Uhlenbeck type. We have
discussed existence of mild solutions and existence, uniqueness, and
ergodicity of invariant measures, generalizing previous work of Vargiolu
\cite{vargiolu}, who considered the situation of a driving Brownian
motion and used a different state space. The choice of the state space
$L^1_{2,\alpha}$ seems to be the standard by now, since its elements
enjoy most desirable features for a forward curve. If the weight
function $\alpha$ grows exponentially at infinity, the HJM dynamics
admits a unique invariant measure. A similar conclusion was obtained
by Tehranchi in \cite{tehranchi}, where a Musiela equation with
state-dependent volatility and Brownian noise was considered, and the
results of \cite{DZ96} could be applied. It would be natural to
consider also in the setting of L\'evy noise a state-dependent
volatility operator, but unfortunately there are comparatively very
few results on evolution equations with jump noise, which need to be
established first.

\subsection*{Acknowledgments}
The financial support of SFB 611 (Bonn) is gratefully acknowledged.
Preliminary results contained in this paper were presented at the
Workshop on Potential Theory and Stochastic Analysis, Bielefeld,
August 2004, and at the Joint Meeting of the Bernoulli Society and the
Institute of Mathematical Statistics, Barcelona, July 2004.

\bibliographystyle{amsplain}

\let\oldbibliography\thebibliography
\renewcommand{\thebibliography}[1]{%
  \oldbibliography{#1}%
  \setlength{\itemsep}{0pt}%
}

\bibliography{ref,merda}

\end{document}